\documentclass[brochure,12pt,french]{smfbourbaki}
        
\usepackage[T1]{fontenc}
\usepackage{lmodern,amssymb,bm}

\usepackage{babel}

\usepackage[ansinew]{inputenc}

\usepackage[colorlinks=true, linkcolor=blue, citecolor=red, urlcolor=blue]{hyperref}

\usepackage[stretch=10,shrink=10,step=2,kerning=true,protrusion=true,expansion=true,final]{microtype}

\usepackage[
backend=bibtex,
style=authoryear, 
citestyle=authoryear-comp,
maxnames=7,
sortcites=false 
]{biblatex}
\usepackage{csquotes}

\newtoggle{tnbcbx@howcited}

\newbibmacro{howcited}{%
  \iftoggle{tnbcbx@howcited}
    {\iffieldundef{shorthand}
       {}
       {\setunit{\addspace}%
        \printtext[parens]{%
          \bibstring{citedas}%
          \setunit{\addcolon\space}%
          \printfield{shorthand}%
          \setunit{\addslash}%
          }}}
    {}}

\renewbibmacro{finentry}{\usebibmacro{howcited}\finentry}

\DefineBibliographyStrings{french}{citedas    = {cit\'e ci-dessus avec l'acronyme}}
\DefineBibliographyStrings{english}{citedas    = {heretofore cited as}}

\DefineBibliographyExtras{french}{\restorecommand\mkbibnamefamily}

\DeclareNameAlias{sortname}{given-family}

\renewcommand{\bibnamedash}{\leavevmode\raise3pt\hbox to3em{\hrulefill}\space}

\AtEveryBibitem{%
  \clearfield{issn} 
  \clearfield{isbn} 
  \clearfield{doi} 
  \clearlist{language} 
  \ifentrytype{online}{}{
  \ifentrytype{unpublished}{}{
    \clearfield{url}
  }
  }
}

\renewbibmacro{in:}{%
    \ifentrytype{article}{}{\printtext{\bibstring{in}\intitlepunct}}}

\DeclareFieldFormat[article,periodical,inreference]{number}{\mkbibparens{#1}}
\DeclareFieldFormat[article,periodical,inreference]{volume}{\mkbibbold{#1}}
\renewbibmacro*{volume+number+eid}{%
    \printfield{volume}%
    \setunit*{\addthinspace}
    \printfield{number}%
    \setunit{\addcomma\space}%
    \printfield{eid}}

\DeclareFieldFormat[article,inbook,incollection]{title}{\enquote{#1}\addcomma} 

\newcommand{\eqdef}{\mathrel{\mathop:}=}

\date{Janvier 2025}
\bbkannee{77\textsuperscript{e} ann\'ee, 2024--2025}  
\bbknumero{1232}                                      

\title{Théorie de l'homotopie quantitative}
\subtitle{d'après Guth, Manin, Weinberger,...}

\author{Pierre Pansu}
\address{Universit\'e Paris-Saclay, CNRS,\\ Laboratoire de math\'ematiques d'Orsay,\\ 91405, Orsay, France.\\ Datashape, Inria Saclay, France.}
\email{pierre.pansu@universite-paris-saclay.fr}

\def\R{\mathbb{R}}
\def\Z{\mathbb{Z}}
\def\Q{\mathbb{Q}}
\def\C{\mathbb{C}}

\def\degre{\mathrm{deg}}
\def\lip{\mathrm{Lip}}
\def\nu#1{\|#1\|_{Lip}}

\usepackage{graphicx}

\begin{document}

\maketitle

\section{Introduction}

Le but de la théorie de l'homotopie, en topologie, c'est de simplifier, après déformation continue, des applications continues entre espaces topologiques. Dans ce texte, on s'intéresse aux questions quantitatives que cela soulève.

Commençons par des exemples. 

\subsection{Nombre de tours}

On s'intéresse aux courbes fermées, i.e. aux applications continues du cercle unité $S^1$ dans le plan privé d'un point $p$ (Figure~1.a). Quand peut-on déformer deux courbes l'une en l'autre en évitant $p$ ? 
\begin{center}
\includegraphics[width=1.8in]{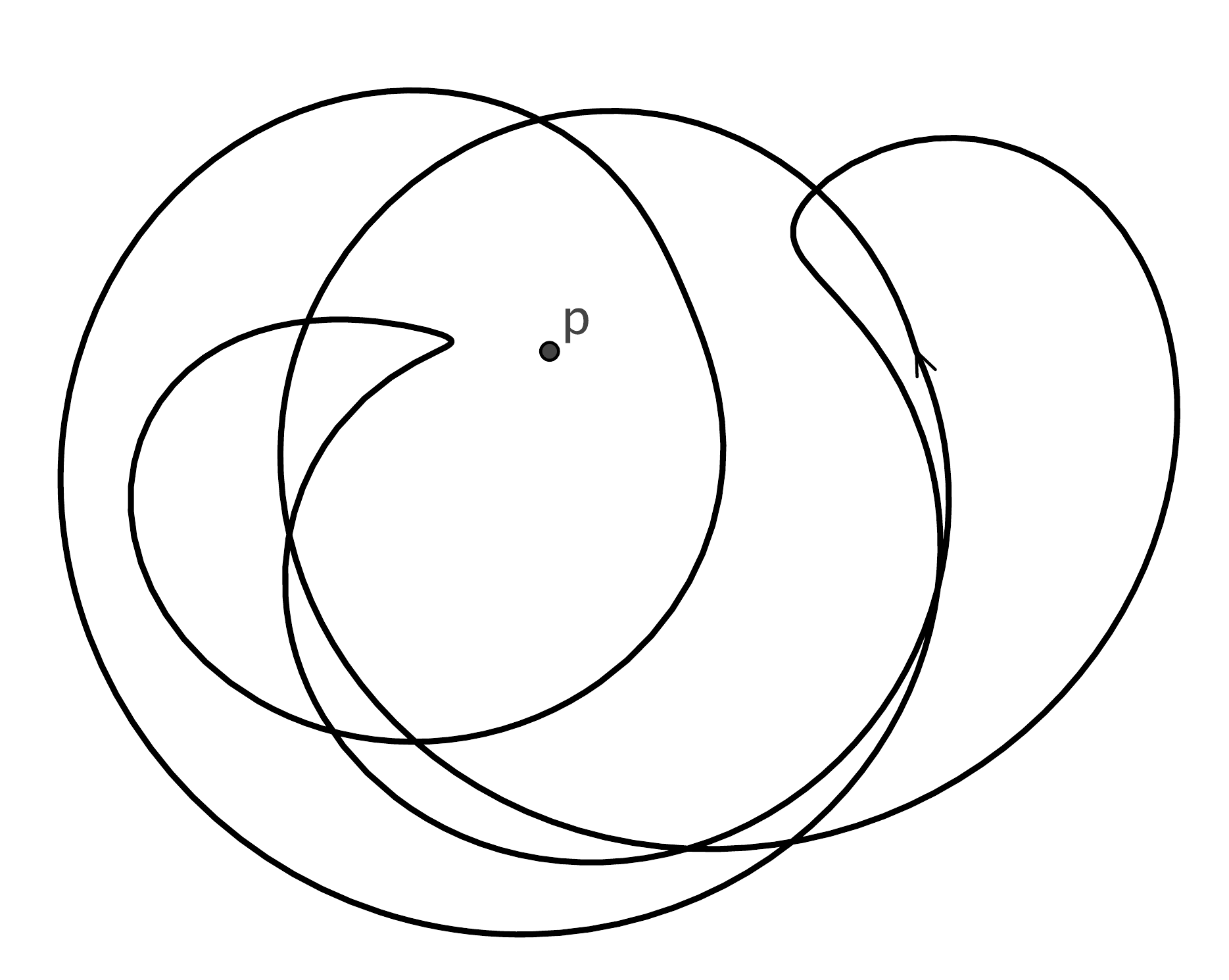}\hskip1cm \includegraphics[width=1.5in]{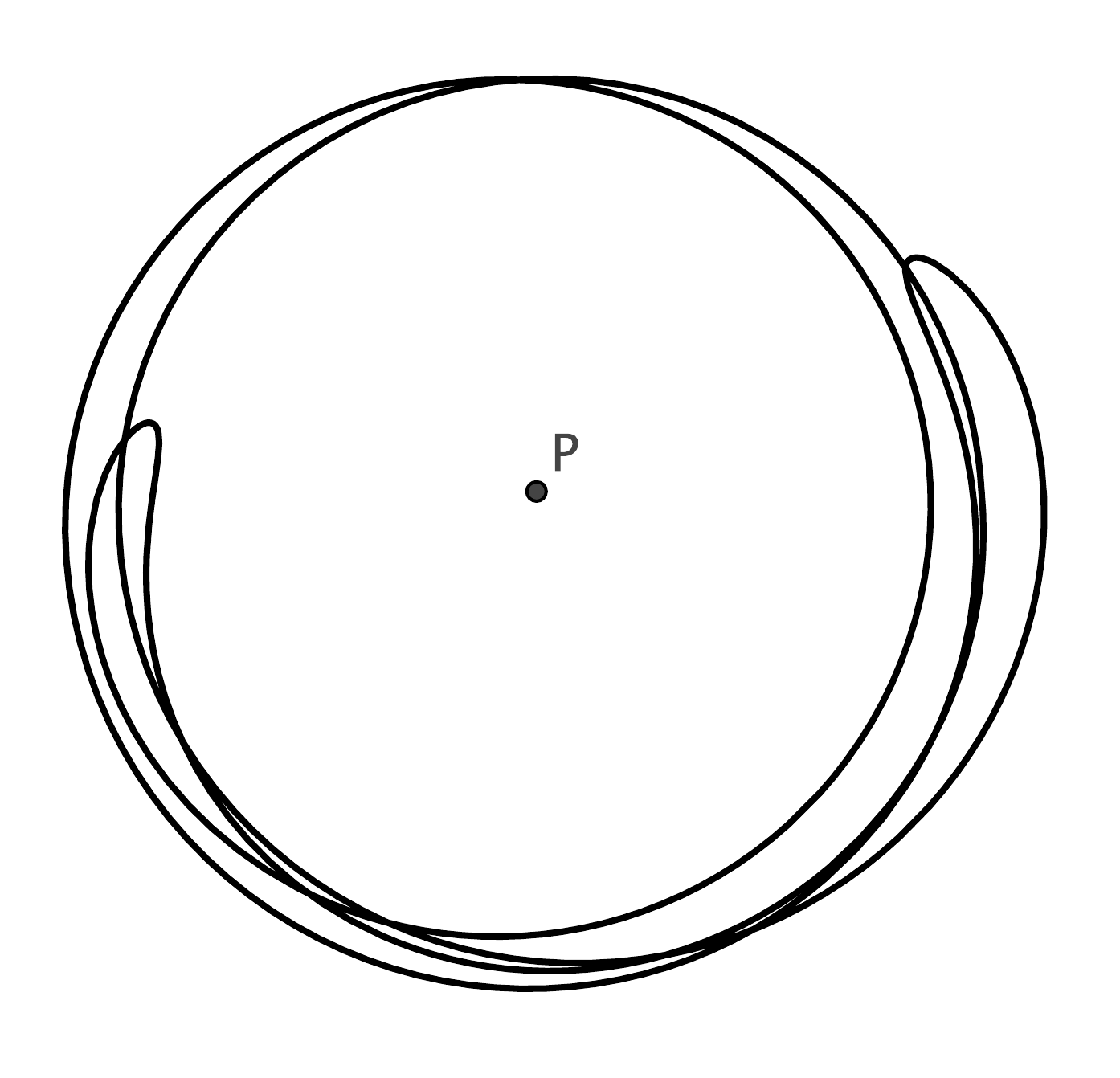}\hskip1cm \includegraphics[width=2in]{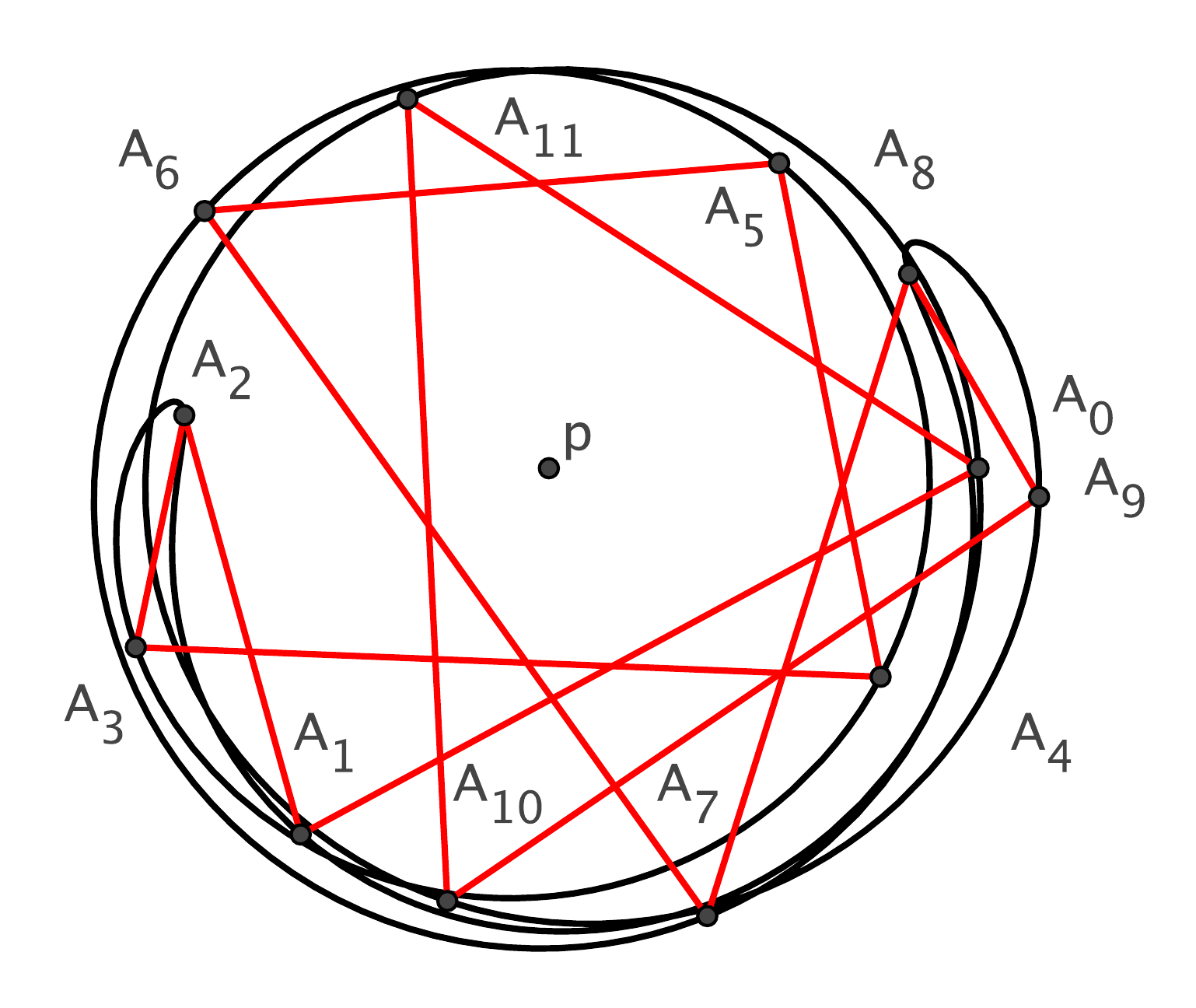}

Figure 1.a \hskip3.3cm Figure 1.b \hskip3.5cm Figure 1.c
\end{center}
Réponse : quand elles font le même nombre de tours autour de $p$. 

Ce nombre est un entier relatif. Est-il compliqué à calculer ? Pour cela, on projette radialement la courbe sur le cercle $\mathcal{C}$ centré en $p$, de rayon~$1$, comme  sur la Figure~1.b. (attention, sur cette figure, pour la lisibilité, on représente une courbe voisine de~$\mathcal{C}$ plutôt que posée sur~$\mathcal{C}$). On va supposer la vitesse angulaire bornée par un nombre entier~$L$, et contrôler la complexité du calcul en fonction de~$L$. On subdivise~$S^1$ en $3L$ intervalles égaux, et on remplace la courbe projetée par une ligne brisée inscrite dans $\mathcal{C}$ (Figure~1.c). 
\begin{center}
\includegraphics[width=1.9in]{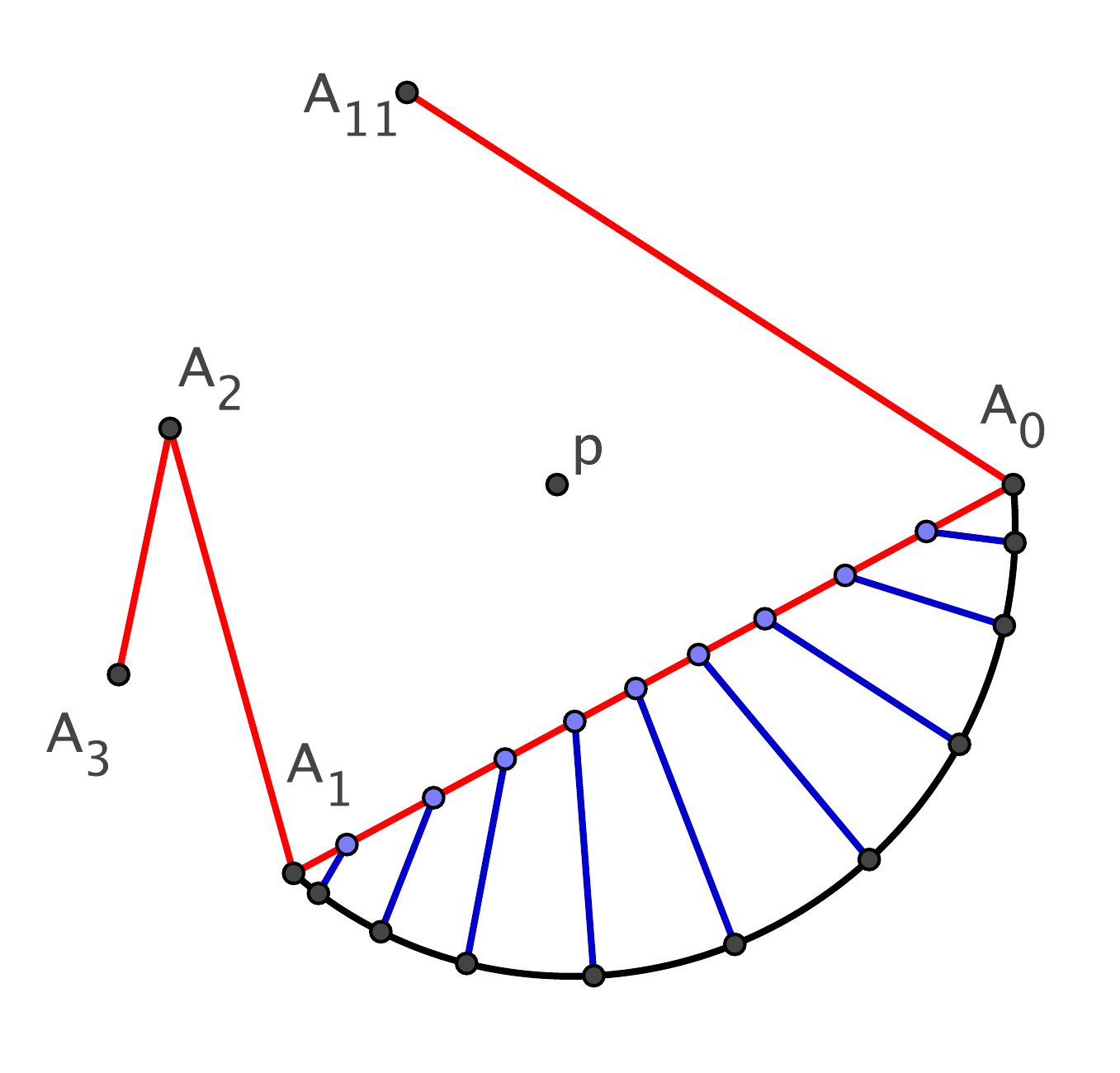}\hskip.5cm \includegraphics[width=1.9in]{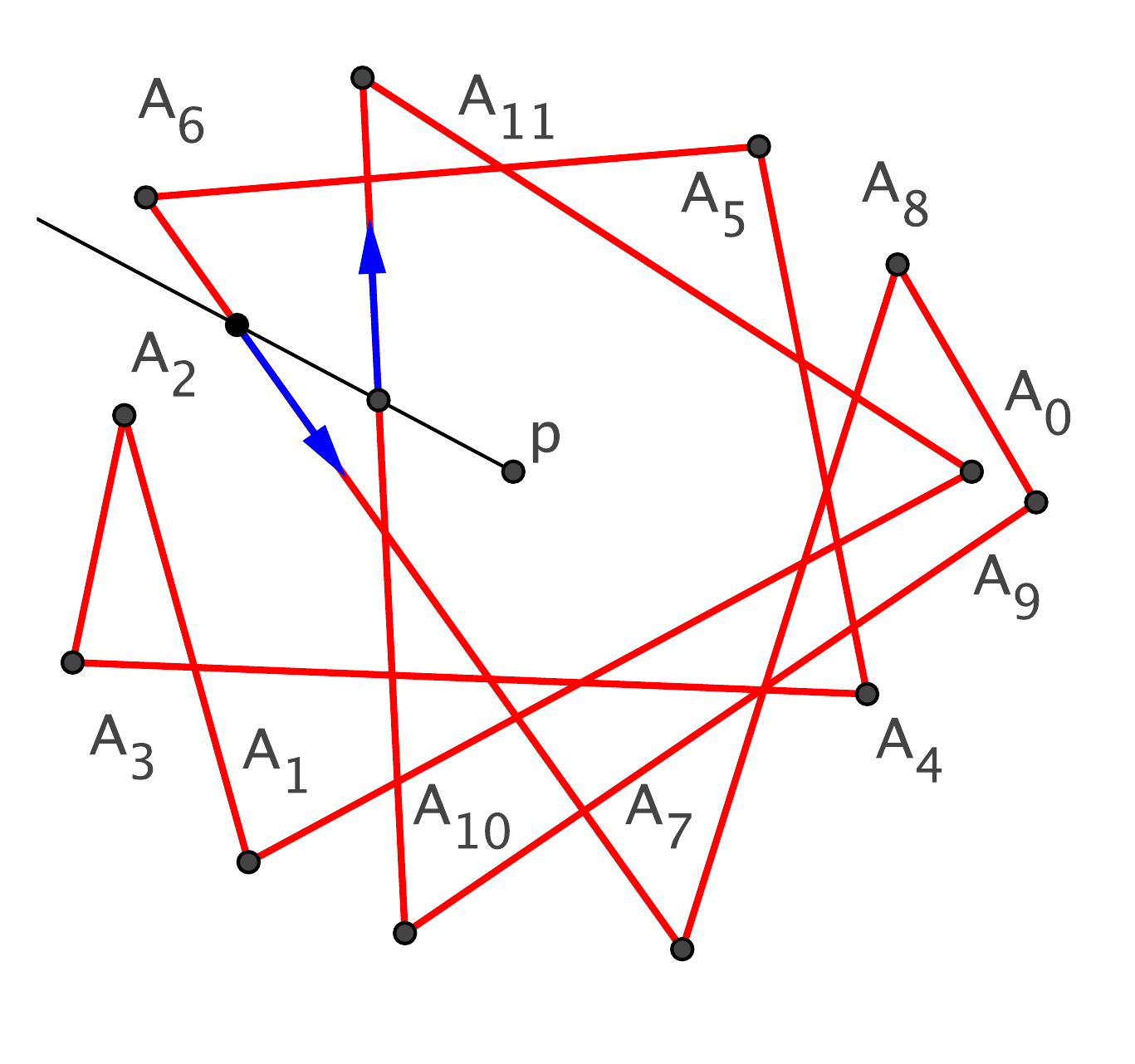}\hskip.5cm \includegraphics[width=1.9in]{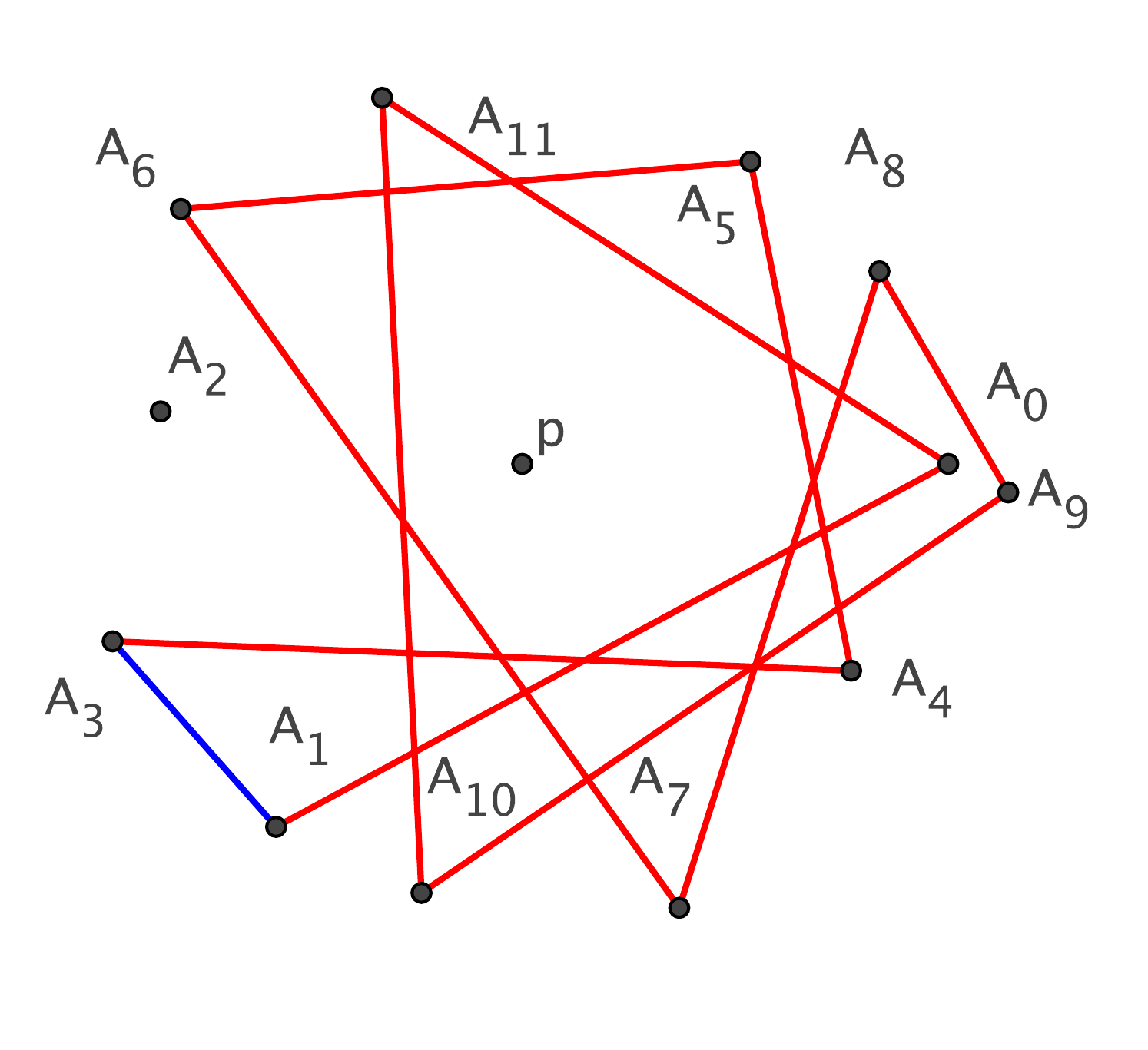}

Figure 2.a \hskip3.5cm Figure 2.b \hskip3.5cm Figure 2.c
\end{center}
Le facteur $3$ fait que chaque segment est vu de $p$ sous un angle au plus $2\pi/3$, on peut donc pousser chaque point de l'arc de $\mathcal{C}$ de mêmes extrémités jusqu'à la corde sans passer par $p$ (Figure 2.a). Cette déformation a lieu à vitesse bornée, prend un temps borné et double au pire la vitesse angulaire. On compte le nombre d'intersections (avec signe) de la ligne brisée avec un rayon quelconque, cela donne le nombre de tours (Figure 2.b). Ici, il vaut $0$. Ce calcul fait intervenir au plus $3L$ fois le même prédicat : intersection d'un segment et d'une demi-droite, son coût est linéaire en $L$.

Pour simplifier la ligne brisée, chaque fois qu'au voisinage d'un sommet $A_i$, l'angle polaire change de sens de variation, on le supprime et on remplace les segments $[A_{i-1}A_{i}]$ et $[A_{i}A_{i+1}]$ par $[A_{i-1}A_{i+1}]$ (Figure 2.c). En moins de $3L$ étapes, on arrive ou bien à un unique sommet (c'est le cas ici), ou bien à une ligne brisée qui coupe tous les rayons dans le même sens, qu'on peut aisément déformer en un paramétrage à vitesse constante de $\mathcal{C}$ parcouru plusieurs fois.

On a vu que la borne supérieure de la vitesse angulaire $L$ constituait une bonne mesure de la complexité de la courbe donnée. Inversement, quels nombres de tours peuvent-ils être réalisés par des courbes de vitesse angulaire $\le L$ ?

Réponse : exactement tous les entiers compris entre $-L$ et $L$. 

Les exemples réalisant les nombres de tours prescrits sont les paramétrages à vitesse constante de $\mathcal{C}$. 

Inversement, le calcul ci-dessus donne une borne $3L$ pour le nombre de tours. La borne optimale $L$ repose sur la géométrie différentielle. Soit $d\theta$ la $1$-forme différentielle angulaire en coordonnées polaires de centre $P$. Si $f\colon S^1\to\mathcal{C}$ désigne la courbe projetée, la forme différentielle (mesurable) $f^*d\theta$ a une norme $\le L$ en chaque point, donc son intégrale $\int_{S^1}f^*d\theta\le 2\pi L$. Or le nombre de tours peut s'écrire $\frac{1}{2\pi}\int_{S^1}f^*d\theta$ (c'est le point de vue adopté par \textcite{Milnor}).

\subsection{Degré}

Le nombre de tours d'une application continue $S^1\to S^1$ se généralise. Soient~$X$ et~$Y$ deux variétés compactes sans bord orientées de même dimension $n$. Le \emph{degré} d'une application lipschitzienne $f\colon X\to Y$ est le nombre entier
$$
\degre(f)\eqdef \frac{\int_{X}f^*\omega}{\int_{Y}\omega},
$$
indépendant du choix de $n$-forme différentielle $\omega$ (L.~\textcite{Brouwer} le définissait différemment). C'est un invariant d'homotopie. Lorsque $Y=S^n$ est la $n$-sphère, c'est le seul invariant: deux applications $X\to S^n$ sont homotopes si et seulement si elles ont même degré (\textcite{Brouwer, Hopf1927} si $X=S^n$). On munit une fois pour toutes~$X$ et~$Y$ de métriques riemanniennes, et on choisit de mesurer la complexité des applications de~$X$ dans~$Y$ par la constante de Lipschitz notée $\lip$. Elle majore le degré,
$$
|\degre(f)|\le C\,\lip(f)^n,
$$
où la constante $C$ ne dépend que des métriques riemanniennes choisies. Lorsque $Y=S^n$, cette borne est optimale : le degré $d$ est réalisable par une application de constante de Lipschitz $O(|d|^{1/n})$. Il suffit pour cela de placer dans $Y$ un nombre $|d|$ boules disjointes de rayon $|d|^{-1/n}$ et d'enrouler chacune sur la sphère, le bord étant envoyé sur un point fixe $o$ de $S^n$.

On verra plus loin (Théorème \ref{th:BGMtri}) que ce n'est pas le cas en général, même lorsque $X=Y$ : il arrive qu'il existe des applications $X\to X$ de degré $d$ élevé, mais leur complexité, mesurée par la constante de Lipschitz, est nécessairement strictement plus grande que $d^{1/n}$.

\subsection{Lacets}

Voici une autre généralisation du nombre de tours. Soit $X$ une variété riemannienne, et $o$~un point de~$X$. Le \emph{groupe fondamental} $\pi_1(X,o)$ est l'ensemble des classes d'homotopie de lacets basés en~$o$, i.e. d'applications de~$S^1$ dans~$x$ astreintes à envoyer le point~$1$ de~$S^1$ sur~$o$. Il possède une structure de groupe induite par la concaténation des lacets. Dans le cas du plan privé d'un point, ce groupe est isomorphe à~$\Z$, l'isomorphisme est donné par le nombre de tours.

Lorsque $X$ est compacte, on peut décrire ce groupe par une présentation finie : un ensemble fini $S$ de générateurs (tout élément du groupe est un produit de puissances de générateurs) et un ensemble fini $R$ de relateurs (tout produit de puissances des générateurs qui vaut l'élément neutre peut se réécrire comme un produit de conjugués de relateurs).

Quand un lacet donne l'élément neutre du groupe fondamental de $X$, on peut le déformer continûment jusqu'au lacet constant, noté encore $o$. 

Est-ce coûteux ? Mesurons le coût en termes de constante de Lipschitz $\lip(H)$ d'une application $H\colon [0,1]\times S^1$ dans $X$ telle que $H_{|\{0\}\times S^1}$ est le lacet $f\colon S^1\to X$ donné, et $H_{|\{1\}\times S^1}$ est constante. La meilleure déformation est celle qui réalise
$$
h(f)\eqdef \inf\{\lip(H)\,;\,H_{|\{0\}\times S^1}=f,\,H_{|\{0\}\times S^1}=o\}.
$$
On s'intéresse à la fonction 
$$
c(\ell)=\sup\{h(f)\,;\,\text{$f$ lacet de longueur }\le \ell, \text{ homotope à }o\}.
$$

On montre aisément que $c$ contrôle la difficulté du \emph{problème du mot} dans la présentation $\langle S|R \rangle$ : à des constantes près, le nombre de conjugués de relateurs nécessaire pour attester qu'un élément $\prod_{i} s_i^{m_i}$ est neutre est au plus $c(\sum_{i}|m_i|)$. Or on sait depuis \textcite{Boone,NovikovP} qu'il existe des présentations de groupes pour lesquelles il n'existe aucun algorithme permettant de décider si un produit de puissances des générateurs est neutre ou non. De telles présentations peuvent être codées dans une variété compacte sans bord de dimension $4$ ou plus. Pour une telle variété, la fonction $c$ n'est pas bornée par une fonction récursive. Familièrement dit, pas vraiment calculable. Ce n'est pas tout de dire qu'un lacet dans une variété est homotope à $0$, réaliser concrètement la déformation continue est une tâche insurmontable dans certaines variétés, dès la dimension $4$.

\subsection{Questions}

Ce qui empêche de simplifier une application continue par déformation, ce sont des invariants homotopiques. Dans ce texte, on va effleurer les questions suivantes :  

- Le calcul des invariants est-il possible (décidable) ? Si oui, à quel coût ?

- Construire des représentants de faible complexité et dont les valeurs des invariants sont prescrites est-il possible ? Si oui, à quel coût ?

- Étant données deux applications qu'on sait pouvoir déformer l'un dans l'autre, quelle est la complexité des déformations nécessaires ?

Les réponses, souvent récentes, sont d'une grande diversité. En outre, bien des questions restent ouvertes, montrant que la topologie n'a pas dit son dernier mot, même en basses dimensions.

\subsection{Remerciements}

Un grand merci à Fedor Manin et Clément Maria pour leur aide, à Larry Guth
pour ses textes très éclairants et à Nicolas Bourbaki pour sa vigilance.

\section{Calculer des invariants}

Si $X$ et $Y$ sont des espaces topologiques, on note $[X,Y]$ l'ensemble des classes d'homotopie d'applications de $X$ dans $Y$. Dans quel sens cet ensemble peut-il être calculé ? Ce paragraphe a pour seule ambition de fournir une sorte d'inventaire journalistique.

\subsection{Déterminer les classes d'homotopie}

Lorsque $X$ est la $n$-sphère $S^n$, $n\ge 2$, $\pi_n(Y)\eqdef [S^n,Y]$ possède une structure de groupe abélien, on l'appelle le $n$-ème groupe d'homotopie de $Y$. Le calcul de tous les groupes d'homotopie rationnelle (i.e. tensorisés par $\Q$), en toute généralité, est $\#$P-difficile \parencite{Anick} (qui, toutefois, suppose $Y$ donné comme un CW-complexe, structure de donnée plus concise qu'une triangulation, par exemple).  La classe $\#$P est l'analogue de NP pour la classe des problèmes de comptage : compter le nombre de solutions d'une équation. Être $\#$P-difficile est considéré comme un résultat strictement plus fort qu'être NP-difficile.

Le calcul de $[X,Y]$ quand $Y$ est $(d-1)$-connexe est décidable si $\mathrm{dim}(X) \le 2d-2$, \parencite{CKMSVW}. Sous ces hypothèses, $[X,Y]$ est un groupe abélien de type fini, on le décrit comme somme directe de groupes monogènes. L'histoire commence avec \textcite{Brown}, qui donne un algorithme pour le calcul des groupes d'homotopie des complexes simpliciaux dont tous les groupes d'homotopie sont finis. La méthode consiste à construire effectivement la tour de Postnikov de $Y$ comme complexe semi-simplicial. La tour est un moyen de décrire un type d'homotopie à partir des briques élémentaires que sont les espaces d'Eilenberg--McLane, espaces dont un seul groupe d'homotopie est non trivial. Cela ramène le calcul de chaque groupe d'homotopie à un calcul de groupe d'homologie, i.e. à des systèmes linéaires sur $\Z$.

La méthode fait intervenir des complexes infinis, car les modèles semi-simpliciaux complets des espaces d'Eilenberg-McLane sont infini (même $K(\Z,1)$ qui a le type d'homotopie du cercle). D'où la nécessité d'introduire une nouvelle catégorie, celle des \emph{espaces à homotopie effective}, structure de donnée implicite qui permet néanmoins des calculs explicites. Le logiciel \emph{Kenzo} implémente ces idées de façon efficace, \parencite{SergeraertGenova}. 

Sous les mêmes hypothèses, il est vraisemblable que le calcul de $[X,Y]$ soit faisable en temps polynômial à $d$ fixé, voir \textcite{CKMVWpol} dûment complété par le caveat de F.~\textcite{Sergeraert}. 

F.~\textcite{Mc} remplace l'hypothèse $d-1$-connexe par l'hypothèse plus vaste ``ayant le type d'homotopie rationnelle d'un H-espace jusqu'en dimension $d$''. Il traite aussi le problème relatif : extension continue à $X$ (de dimension $\le 2d-1$) d'une application donnée sur un sous-complexe $A$ de $X$.

Dès qu'on quitte ce régime ``stable'', on s'attend à ce que le calcul de $[X,Y]$ devienne indécidable, et c'est effectivement le cas pour le problème relatif, \parencite{CKMVWund}: \emph{il existe un espace $d-1$-connexe $Y$ qui a la propriété suivante. Le problème de l'existence d'un prolongement de $A$ à $X$, pour tout complexe simplicial $X$ et $A\subset X$ de dimension $2d$ et toute application continue $f\colon A\to Y$, est indécidable. }La méthode, qui s'appuie sur la solution du 10ème problème de Hilbert, consiste à coder un système d'équations diophantiennes quadratiques dans un problème d'extension. C'est possible à partir de la dimension $2d$, $d$ pair, parce que les sphères paires ont deux groupes d'homotopie infinis, $\pi_d(S^d)$ et $\pi_{2d-1}(S^d)$, avec une multiplication commutative de l'un dans l'autre. $Y=S^d$ fait l'affaire dans ce cas. Lorsque $d$ est impair, on remplace $S^d$ par $Y=S^d \vee S^d$, et un produit de Whitehead de $\pi_{d}(S^d \vee S^d)$ dans $\pi_{2d-1}(S^d \vee S^d)$. Ce dernier est alterné, mais tout système d'équations diophantiennes peut se ramener aussi à un système bilinéaire antisymétrique.

\subsection{Cas de la dimension $3$}

En dimension $3$, on dispose d'une floraison d'invariants, découlant des multiples façons de construire des $3$-variétés à partir de briques élémentaires.

Un \emph{corps en anses de genre $g$} est une variété à bord délimitée par une surface de genre $g$ plongée de façon usuelle dans $\R^3$. En recollant deux corps à anses au moyen d'un difféomorphisme du bord renversant l'orientation, on obtient une $3$-variété, munie d'une \emph{décomposition de Heegard}. Le plus petit genre d'une décomposition de Heegard s'appelle le \emph{genre de Heegard} de la $3$-variété. Le calcul du genre de Heegard est décidable mais NP-difficile \parencite{BDTS}. La classe NP est la classe des problèmes de décision pour lesquels on peut fournir un certificat vérifiable en temps polynômial. Étant donné un entier $g$, c'est l'assertion ``genre de Heegard $\le g$'' à laquelle on peut réduire en temps polynômial tout problème de la classe NP. La démonstration consiste à coder des formules booléennes en forme conjonctive (i.e. avec des ``et'', des ``ou'' et des négations) dans des décompositions de Heegard de $3$-variétés.

Le genre d'un noeud dans une $3$-variété est le plus petit genre d'une surface orientable plongée dont le bord est le noeud donné. Il vaut $0$ exactement quand le noeud n'est pas noué. Le calcul du genre d'un noeud dans une $3$-variété quelconque est NP-complet. Dans une 3-variété fixée, c'est coNP \parencite{AHT, LY}. Autrement dit, ``genre du noeud $\le g$'' est dans NP (et même NP-complet), et le contraire, ``genre du noeud $> g$'' est aussi dans NP si on se limite à une $3$-variété fixée. En particulier, pour les noeuds dans la sphère, décider qu'un noeud est noué est dans~NP.

Le calcul des invariants quantiques des $3$-variétés jette un pont entre topologie et algorithmique quantique. Les invariants $\tau$ dont il est ici question obéïssent à une règle de fusion : un espace vectoriel $\tau(S)$ est associé à toute surface orientée $S$,  de sorte que $\tau(-S)=\tau(S)^*$, un vecteur $\tau(V)\in\tau(\partial V)$ pour toute $3$-variété à bord, de sorte que pour une $3$-variété fermée $X$, divisée en deux moitiés $V_1$ et $V_2$ par une surface $S$, l'invariant numérique $\tau(X)$ soit donné par $\tau(X)=\langle \tau(V_1),\tau(V_2) \rangle$. Pour les noeuds et entrelacs, les surfaces sont des sphères avec points marqués, les moitiés sont des boules munies d'arcs reliant les points marqués du bord (un peu comme sur la Figure 3). Le prototype est le polynôme de Jones $V(t)$, $t\in\C$, $|t|=1$. Dans ce cas, l'espace vectoriel $V(2n)=V($sphère marquée de $2n$ points$)$ est de dimension finie. Il faut éventuellement passer à un quotient lorsque $t=e^{2i\pi/r}$ est une racine de l'unité. Le groupe des tresses $B_{n}$ agit par matrices unitaires sur $V(2n)$. Lorsque l'image de cette représentation unitaire est un sous-groupe dense du groupe unitaire, tout circuit quantique peut être simulé approximativement par une tresse, la probabilité que le circuit réponde oui est proche du nombre $|V(t)|^2$, au sens où le rapport des deux est proche de $1$, \parencite{FKW}. Cela permet de montrer que le calcul approché ---au sens multiplicatif--- de $|V(e^{2i\pi/r})|^2$ pour un noeud, si $r=5$ ou $r\ge 7$, est postBQP-difficile, \parencite{Kuperberg, AL}, ce qui, en vertu d'un théorème général de théorie de la complexité \parencite{Aaronson}, coïncide avec $\#$P-difficile. On pense que le calcul reste difficile même si on limite la classe de noeuds en jeu, par exemple en bornant certains invariants du noeud. La méthode s'étend à certains invariants quantiques de $3$-variétés. En revanche, si $r=1,2,3,4,6$, l'image de la représentation de Jones du groupe des tresses n'est pas dense. Dans cette situation de non densité, le calcul exact d'un invariant quantique peut être $\#$P-difficile (voir \textcite{KiMe} pour l'invariant de Witten-Reshetikhin-Turaev  en $e^{2i\pi/4}$ pour le groupe $SU(2)$), mais on pense qu'un calcul approché est moins difficile. En outre, des algorithmes polynômiaux existent pour certains invariants de cette famille si on fixe $b_1$ \parencite{MS, DMS}. On peut se rabattre sur une approximation additive : pour tout $\epsilon>0$, l'algorithme fournit avec forte probabilité une approximation $\hat\tau$ de $\tau$ telle que $|\hat\tau-\tau|<\epsilon$. Cela ne permet pas de décider si $\tau=0$ ou non. De tels algorithmes polynômiaux existent à condition d'utiliser de l'aléa quantique, le problème de l'approximation additive appartient donc à la classe BQP \parencite{FLW1, FLW2}. 

Des travaux de W. Thurston, il résulte que la plupart des noeuds premiers de la $3$-sphère ont un complémentaire difféomorphe à une (unique) variété hyperbolique de volume fini. Ce volume est un invariant important. Une façon de le trouver consiste à construire une triangulation idéale (les sommets sont rejetés à l'infini) du complémentaire du noeud, à réaliser chaque tétraèdre dans l'espace hyperbolique (i.e. choisir un birapport, un nombre complexe) et à résoudre ---de façon approchée--- les équations algébriques dictées par le recollement. Cette heuristique a été implémentée par \textcite{SnapPy} dans le logiciel \emph{SnapPy}, avec succès (le logiciel trouve toujours une solution après quelques essais, voir l'étude expérimentale d'O.~\textcite{Owen}), sans qu'on puisse le démontrer rigoureusement. 

\subsection{Reconnaître la trivialité des invariants}

À défaut de calculer des invariants, on peut être moins ambitieux et se contenter de savoir s'ils sont triviaux ou non.

La trivialité du groupe fondamental est indécidable pour les variétés de
dimension $\ge 5$, 
\parencite{Adyan, Rabin}. Par conséquent, reconnaître la $n$-sphère, $n\ge 5$, est indécidable, \parencite{NovikovP}.

Reconnaître la 3-sphère, en vertu de la conjecture de Poincaré, équivaut à reconnaître la trivialité du groupe fondamental. Ce problème est non seulement décidable, il est dans NP \parencite{Rubinstein, Thompson, Schleimer}. L'algorithme repose sur les notions de surface normale et de surface presque normale dans une $3$-variété triangulée. Une surface normale est l'analogue combinatoire d'une surface minimale stable en géométrie riemannienne. Une surface presque normale est l'analogue combinatoire d'une surface minimale instable. Reconnaître la trivialité de l'homologie est dans P, donc on peut se focaliser sur les $3$-variétés dont l'homologie est nulle (les ``sphères d'homologie''). Pour celles-ci, le certificat est une surface presque normale telle que le bord d'aucune des composantes du complémentaire n'est une sphère normale, sauf si c'est le voisinage d'un sommet. Les surfaces normales et presque normales peuvent être codées par des vecteurs entiers sujets à des équations, on peut vérifier aisément leur nature, d'où le caractère NP. On peut les énumérer, d'où la décidabilité. 

Reconnaître la 3-sphère est dans coNP, conditionnellement à l'hypothèse de Riemann généralisée (GRH), \parencite{Zentner}. La méthode est totalement différente. Pour les sphères d'homologie, le certificat est une représentation non triviale du groupe fondamental dans $SL(2,\mathbb{C})$. Or $Hom(\pi_1,SL(2,\mathbb{C}))$ est une variété algébrique définie sur $\Z$. Le rôle de GRH est de trouver un point de cette variété défini sur $\Z/p\Z$ pour un nombre premier $p$ pas trop grand (voir \textcite{KuperbergNP}). 

La combinaison des deux résultats, NP et coNP, signifie probablement que reconnaître la $3$-sphère n'est pas NP-complet, et relève donc d'un degré de complexité intermédiaire.

L'énumération ci-dessus peut donner l'impression que l'indécidabilité ne se produit pas en topologie en dimensions $\le 3$. Ce n'est pas le cas, puisque reconnaître les groupes fondamentaux de $3$-variétés parmi les présentations finies de groupes est indécidable, \parencite{Stallings}. Le problème de décider si les groupes fondamentaux de deux $3$-variétés sont isomorphes est encore partiellement ouvert, \parencite{AFW}.

Reconnaître la $4$-sphère est-il décidable ? La question est ouverte à ce jour. La méthode employée en dimension $3$ ne s'étend pas : il existe des sphères d'homologie de dimension $4$ distinctes de $S^4$ sans représentations dans $SL(2,\mathbb{C})$, \parencite{Zentner}. La méthode employée en dimension $\ge 5$ donne davantage d'espoir. Un résultat de M.~\textcite{Kervaire} affirme que tout groupe parfait (d'abélianisé trivial) de présentation équilibrée (même nombre de générateurs et de relateurs) est le groupe fondamental d'une sphère d'homologie de dimension $4$. Il suffirait de disposer de familles de présentations équilibrées de groupes parfaits dont la trivialité est indécidable, mais on n'en connaît pas. En revanche, il y a des familles de présentations équilibrées du groupe trivial pour lesquelles le nombre de conjugués de relateurs à multiplier pour montrer que les générateurs sont triviaux est très grand \parencite{Bridson,Lishak}, tout en restant récursif. Cela donne des bornes inférieures sur la taille de l'espace des triangulations de certaines $4$-variétés \parencite{LN}. Noter que décider si deux variétés de dimension $4$ sont homéomorphes est indécidable, \parencite{Markov}.

\section{Norme des classes d'homotopie}

\subsection{Paysage de Morse}

On fixe une fois pour toutes les variétés riemanniennes $X$ et $Y$ et, suivant M.~\textcite{GromovHdil, GromovQH}, on s'intéresse au paysage ``de Morse'' de la fonction $\lip$ sur l'espace $\lip(X,Y)$ des applications lipschitziennes, i.e. à la famille des ensembles de sous-niveaux
$$
\lip_{\le L}\eqdef \{f\in \lip(X,Y)\,;\,\lip(f)\le L\}.
$$
Étant donnée une classe d'homotopie $[f]\in[X,Y]$, on peut penser à la plus petite constante de Lipschitz d'un représentant de la classe, notée $\nu{[f]}$, comme à une sorte de norme.
\begin{itemize}
  \item Combien de composantes connexes de $\lip(X,Y)$ intersectent $\lip_{\le L}$ ? De façon équivalente, quelle est la croissance des boules dans $[X,Y]$ pour la norme $\nu{\cdot}$ ?
  \item Est-ce que deux éléments de $\lip_{\le L}$ qui sont connectés dans $\lip(X,Y)$ le sont déjà dans $\lip_{\le \Lambda(L)}$, pour une fonction $\Lambda$ à trouver ? Par exemple, lorsque $X=Y=S^1$, $\Lambda(L)=L$ convient.
\end{itemize}

\subsection{Un exemple : l'invariant de Hopf}

L'étude de $\pi_3(S^2)$ remonte à \textcite{Hopf1931}. H. Hopf introduit la \emph{fibration de Hopf} $\phi$, qui envoie un vecteur unitaire de $\C^2$ (i.e. un élément de $S^3$) sur la droite complexe qu'il engendre (un élément de la droite projective complexe, difféomorphe à $S^2$). Ses fibres sont des grands cercles de $S^3$ qui sont enlacés : tout disque bordé par une des fibres coupe toutes les autres fibres. Hopf en déduit que $\phi$ n'est pas homotope à $0$. De nos jours, suivant \textcite{Whitehead}, on procèderait comme suit. Soit $f\colon S^3\to S^2$ une application lipschitzienne. Soit $\omega$ une $2$-forme différentielle sur $S^2$ d'intégrale non nulle. La forme $f^*\omega$ est fermée. Comme la $2$-cohomologie réelle de $S^3$ s'annule, il existe une $1$-forme différentielle $\alpha$ telle que $d\alpha=f^*\omega$. Le nombre
$$
h(f)\eqdef \frac{\int_{S^3}\alpha\wedge f^*\omega}{\int_{S^2}\omega}
$$
est un entier appelé \emph{invariant de Hopf}. Il ne dépend que de la classe d'homotopie de~$f$. Il vaut~$1$ pour la fibration de Hopf~$h$ et $0$~pour l'application constante, donc elles ne sont pas homotopes. On montre que c'est le seul invariant d'homotopie, et que c'est un homomorphisme de groupes, donc $\pi_3(S^2)=\Z$. 

Par construction, $\|f^*\omega\|_{\infty}\le \lip(f)^2\|\omega\|_{\infty}$. On montre qu'on peut choisir $\alpha$ telle que $\|\alpha\|_{\infty}\le C\,\|f^*\omega\|_{\infty}$, d'où
$$
|h(f)|\le C\,\lip(f)^4.
$$
Inversement, en composant une application $f$ d'invariant de Hopf $h\in\Z$ avec une application $g\colon S^2\to S^2$ de degré $k$, on réalise $h(g\circ \phi)=k^2 h$ et $\lip(g\circ f)\le \lip(g)\lip(f)\le C\,k^{1/2}\lip(f)$. Autrement dit, pour tous entiers $k>0$ et $h$, $\nu{k^2 h}\le C\,k^{1/2}\nu{h}$. Avec la sous-additivité à constante multiplicative près, cela entraîne que $\nu{h}=O(|h|^{1/4})$. Le nombre de classes d'homotopie contenant un représentant de constante de Lipschitz $\le L$ est donc de l'ordre de $L^4$.

\subsection{Majoration du nombre de classes de complexité donnée}

Dès les années 70, M.~\textcite{GromovHdil} a eu l'intuition que la théorie des modèles minimaux de D.~\textcite{Sullivan}, permettait de généraliser l'exemple de l'invariant de Hopf, au moins dans le cas où $X$ est une sphère. Ce programme n'a été complété que récemment.

\begin{theo}[\cite{MWmc}]\label{th:MWmc}
On suppose $Y$ simplement connexe (i.e. son groupe fondamental n'a qu'un élément). La croissance des boules dans $[X,Y]$ est au plus polynômiale.
\end{theo}
Le théorème est valable plus généralement pour les espaces $Y$ \emph{nilpotents}, i.e. tels que $\pi_1(Y)$ est nilpotent et agit de façon nilpotente sur les groupes d'homotopie supérieurs.

La théorie des modèles minimaux permet de cerner l'homotopie \emph{rationnelle} d'une variété $X$ triangulée (en particulier, ses groupes d'homotopie tensorisés par $\Q$) au moyen d'algèbres différentielles graduées (dga). Par exemple, il existe une version rationnelle, notée $A^*X$, de l'algèbre des formes différentielles. Parmi toutes les dga possibles, il y en a une qui est minimale. Elle est notée $\mathcal{M}^*_X$, et est étroitement reliée à la tour de Postnikov de $X$. Comme cette dernière, le modèle minimal est une tour qui permet des raisonnements par récurrence. 

L'homotopie rationnelle de la variété $Y$ s'incarne dans un espace $Y_\Q$, de sorte que, par exemple, $\pi_k(Y_\Q)=\pi_k(Y)\otimes\Q$. Il y a une bijection entre l'ensemble des classes d'homotopie $[X,Y_\Q]$ et l'ensemble des classes d'homotopie de morphismes de dga entre $\mathcal{M}^*_Y$ et $A^*X$, et ces derniers ont une description en termes finis (bien que la dga $A^*X$ soit de dimension infinie). Comme chaque nouvel étage du modèle minimal $\mathcal{M}^*_Y$ ajoute des relations de la forme $da=b$, l'analyse de la norme d'un morphisme se déroule comme dans l'exemple de l'invariant de Hopf, conduisant à des bornes polynômiales.

Pour tout $X$, il existe une application $[X,Y]\to [X,Y_\Q]$ dont les fibres sont finies, \parencite{Sullivan}, mais de cardinal non borné uniformément en général. \textcite{MWmc} démontrent que le cardinal de la fibre d'un morphisme croît au plus polynômialement en fonction de la norme du morphisme. Cela complète la démonstration du Théorème \ref{th:MWmc}.

Le Théorème \ref{th:MWmc} ne dit pas que le nombre de classes d'homotopie d'applications croît de façon exactement polynômiale. D'ailleurs, ce n'est pas vrai en général. Les auteurs donnent un exemple où la croissance est en $L^8 \log L$. Il s'agit d'un complexe obtenu en ajoutant une ou deux cellules de dimension plus grande à un bouquet de deux sphères, de sorte que les classes d'homotopie puissent être caractérisées par des sortes d'invariants de Hopf relatifs: la restriction de l'application donnée $f$ au bouquet est un couple d'applications $(f_1,f_2)$. On s'arrange pour disposer d'un invariant de Hopf de $f_1$ à valeurs dans $\Z/h(f_2)\Z$, d'où une boule qui fibre sur une honnête boule vectorielle (de cardinal exactement polynômial) avec des fibres de cardinaux variables. En exploitant ce mécanisme, on obtient le théorème suivant.

\begin{theo}[\cite{Mz}]\label{th:Mz}
Pour tout rationnel $r>4$, il existe des espaces simplement connexes $X$ et $Y$ tels que le nombre de classes d'homotopie d'applications de $X$ dans $Y$ telles que $\nu{\cdot}\le L$ croît comme $L^r$ (au sens où le rapport des $\log$ tend vers $r$).
\end{theo}

\subsection{Minoration du nombre de classes de complexité donnée}

Il s'agit de construire des exemples d'applications contournées sur le plan homotopique, mais de constante de Lipschitz modérée. Ce n'est pas simple, même quand il s'agit d'envoyer une variété $X$ de dimension $n$ sur $Y=X$. Dans ce cas, on s'intéresse à la relation entre degré et constante de Lipschitz. L'inégalité 
$$
\degre\le C\,\lip^n
$$
est-elle optimale ? Peut-être pas. Il pourrait y avoir une interférence de l'homologie en degré $k<n$ ; par exemple, que toutes les applications $X\to X$ de degré $d$ multiplient une classe de degré $k$ par un facteur $>d^{k/n}$. 

\begin{defi}[\cite{BM}]\label{def:scalable}
On dit qu'une variété est \emph{dilatable} (en anglais, scalable), si elle admet des applications dans elle-même de degrés $d$ arbitrairement grands et de constantes de Lipschitz $O(d^{1/n})$.
\end{defi}

La dilatabilité a une caractérisation en termes d'homotopie rationnelle et même réelle. Observons d'abord que la dilatabilité entraîne la propriété suivante, qui a fait son apparition en homotopie rationnelle dès les années 1970 (avec une définition initiale différente).

\begin{defi}[\cite{Sullivan}]\label{def:formal}
On dit qu'une variété $X$ est \emph{formelle} (en anglais, formal), si elle admet une application dans elle-même qui induit sur chaque espace d'homologie $H_k(X;\R)$ la multiplication par $p^k$, pour un $p\not=0$.
\end{defi}

\begin{theo}[\cite{BGM}]\label{th:BGMsc}
Une variété $X$ de dimension $n$ est dilatable si et seulement si elle est formelle et l'anneau de cohomologie $H^*(X,\R)$ se plonge dans l'algèbre extérieure $\bigwedge^* \R^n$.
\end{theo}

Les variétés $(n-1)$-connexes de dimension $2n$ (par exemple, les variétés de dimension $4$ simplement connexes) sont dilatables si et seulement si les deux entiers $k$ et $\ell$ constituant la signature $(k,\ell)$ de la forme d'intersection sur $H^n(X,\R)$ sont inférieurs ou égaux à $\frac{1}{2}{2n \choose n}$. 
Quand on fait des sommes connexes, la cohomologie s'additionne, et quand le nombre de termes est assez grand, la variété obtenue n'est plus dilatable.

La dilatabilité permet de construire des homotopies très bien contrôlées.

\begin{theo}[\cite{BM}]\label{th:BMhomot}
Si une variété $Y$ est dilatable, les applications de tout complexe $X$ vers $Y$, de constante de Lipschitz $L$, homotopes à $0$, le sont par des homotopies de constantes de Lipschitz linéaires en $L$.
\end{theo}

Nous avons vu que dilatable implique formel, mais la réciproque n'est pas vraie. Ces deux notions sont les clés de la relation entre degré et constante de Lipschitz, qui prend la forme d'une trichotomie.

\begin{theo}[\cite{BGM}]\label{th:BGMtri}
Soit $X$ une variété de dimension $n$. On s'intéresse au degré maximal $\delta(L)$ d'une application de $X$ dans $X$ de constante de Lipschitz $\le L$.
\begin{enumerate}
  \item Si $X$ est formelle et dilatable, alors pour une infinité de valeurs de $L$,
  $$
C\,L^{n} \le \delta(L) .
  $$
  \item Si $X$ est formelle et non dilatable, alors il existe des constantes $a$ et $b$ ne dépendant que de l'anneau de cohomologie $H^*(X,\R)$ telles que, 
  $$
L^n(\log L)^{-b} \le  \delta(L)\le L^{n}(\log L)^{-a} .
  $$
  \item Si $X$ n'est pas formelle, alors il existe une constante $\alpha<n$ telle que 
  $$
\delta(L)\le L^{\alpha}.
  $$
\end{enumerate}
\end{theo}

Il s'agit d'améliorer l'estimation naïve sur le degré. La démonstration consiste à faire une analyse de Fourier du tiré en arrière de la forme volume. Cette analyse suggère aussi la construction d'exemples. Si les très basses fréquences dominent, cela fournit un plongement de la cohomologie dans l'algèbre extérieure. Sinon, on gagne un facteur polylogarithmique sur l'estimation naïve $\delta(L) \le C\, L^n$. 

Ces méthodes apportent un élément de réponse à une question de M.~\textcite{GromovQH}. Disons qu'une classe d'homotopie $\alpha\in\pi_n(Y)$ est \emph{distordue} si les multiples $k\alpha$ possèdent des représentants de constante de Lipschitz plus faible que les $O(k^{1/n})$ attendus. M. Gromov demande si une classe est distordue si et seulement si la classe d'homologie correspondante dans $H_n(Y,\Z)$ est nulle, et si, dans ce cas, on gagne sur l'exposant, $1/n$ étant remplacé par $1/(n+1)$. F.~\textcite{Mpc} le démontre lorsque $Y$ est un espace symétrique compact, variétés qui sont dilatables. En revanche, la réponse est négative pour $n=5$ et la somme connexe de $4$ copies de $\C P^2 \times S^2$, privée d'un point, variété qui n'est pas dilatable : la constante de Lipschitz ne descend pas jusqu'à $O(k^{1/6})$ \parencite{BM}.

\subsection{Contrôle des homotopies}

On en vient à la question de construire des homotopies contrôlées (en constante de Lipschitz) entre applications lipschitziennes homotopes. L.~\textcite{Guth} souligne la similitude de ce problème avec l'inégalité isopérimétrique : la donnée est sur le bord d'une variété (le produit $[0,1]\times X$) et il s'agit de contrôler ce qui se passe à l'intérieur. L'inégalité isopérimétrique est apparue plus haut, quand il s'agissait d'estimer l'invariant de Hopf. En effet, trouver une primitive bornée d'une forme différentielle fermée bornée, par dualité, équivaut à remplir un cycle en contrôlant le volume du remplissage. Elle joue un rôle dans chacun des résultats de cette section.

Sous l'hypothèse assez faible que l'espace d'arrivée est simplement connexe, le problème de contrôle des homotopies a une solution polynômiale.

\begin{theo}[\cite{Mpc}]\label{th:Mpc}
Lorsque $Y$ est simplement connexe, il existe un polynôme $P_{X,Y}$ tel que si deux éléments de $Lip_{\le L}$ sont homotopes, alors il existe une homotopie entre elles de constante de Lipschitz $\le P_{X,Y}(L)$. Le polynôme $P_{X,Y}$ est même quadratique quand l'une des deux applications est constante. 
\end{theo}
Autrement dit, la fonction $\Lambda$ qui garantit une borne de Lipschitz pour les homotopies peut être choisie polynômiale, ou même quadratique.

La méthode utilise l'homotopie rationnelle. Étant données deux applications homotopes, F. Manin construit d'abord une homotopie au niveau des morphismes de dga dont la norme est contrôlée par $L$. Il la corrige aux deux bouts pour qu'elle reflète une véritable homotopie. Enfin, il construit pour tout morphisme qui provient d'une classe d'homotopie entre espaces un bon représentant dont la constante de Lipschitz est bornée linéairement en fonction de la norme du morphisme.

Le résultat est plus précis : il stipule que les homotopies construites sont de \emph{longueur} bornée, i.e. chaque point de $X$ décrit au cours de la déformation un chemin de longueur uniformément bornée dans $Y$.

La borne polynômiale peut-elle être améliorée ? Certainement pas en toute généralité, mais sous des hypothèses plus restrictives sur $Y$, la réponse est oui. On a déjà rencontré ce phénomène au Théorème \ref{th:BMhomot}. Voici deux autres résultats dans ce sens.

\begin{theo}[\cite{FW}]\label{th:FW}
Lorsque $Y$ a ses groupes d'homotopies finis jusqu'en dimension $d$, la fonction $\Lambda$ peut être choisie linéaire. Plus précisément, il existe une constante $C_{d}$ telle que, pour tout complexe simplicial $X$ de dimension $\le d$, si deux éléments de $Lip_{\le L}(X,Y)$ sont homotopes, alors il existe une homotopie entre elles de constante de Lipschitz $\le C_{d}\,L$. Réciproquement, si de telles homotopies linéaires existent, $Y$ a ses groupes d'homotopies finis jusqu'en dimension $d$
\end{theo}

\begin{theo}[\cite{CDMW}]\label{th:CDMW}
Supposons que $\mathrm{dim}(X)\le n$ et que $Y$ est, jusqu'en dimension $n$, rationnellement équivalent à un produit d'espaces d'Eilenberg--McLane. Alors la fonction $\Lambda$ peut être choisie linéaire.
\end{theo}
On ne peut pas faire mieux ! L.~\textcite{Guth} décrit ainsi les arguments dans le cas particulier où $X=S^n$ et $Y=S^m$. 

Lorsque $n<m$, on montre que toute sphère $L$-lipschitzienne évite une boule pas trop petite, il ne reste qu'à rétracter le complémentaire de cette boule.

Lorsque $n=m$, d'après la définition de Brouwer du degré, étant donnée une application lisse homotope à $0$, il existe une fibre constituée d'un nombre pair de points affectés de signes $+$ et $-$ en quantités égales. Construire une homotopie, c'est apparier les $+$ et les $-$ et connecter les paires dans la boule $B^{n+1}$ ou, ce qui revient au même, dans un cube. 

\begin{center}
\includegraphics[width=4in]{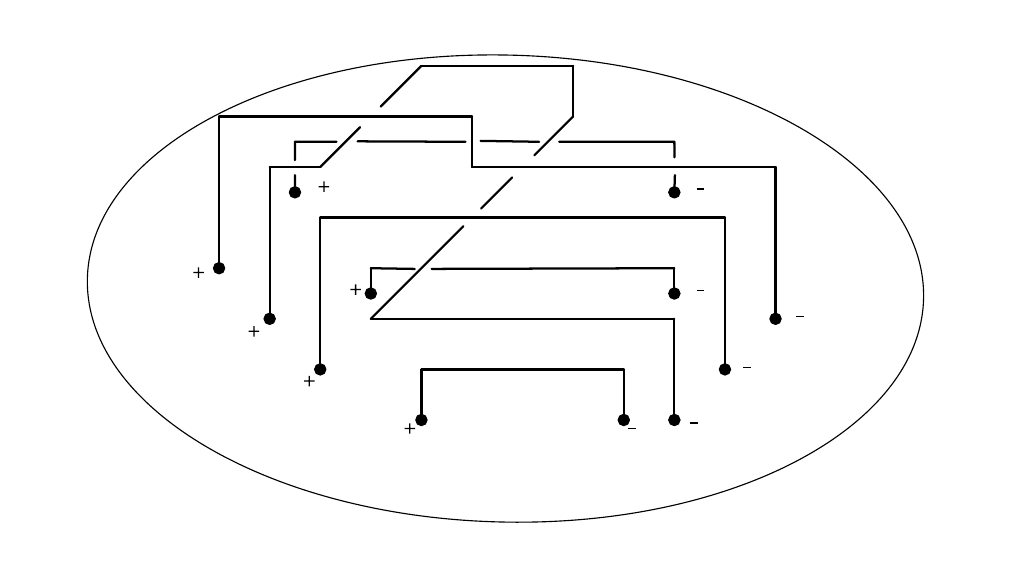}

Figure 3
\end{center}

Construire une homotopie contrôlée, c'est connecter les paires par des chemins éloignés les uns des autres. La construction de tels chemins est due à \textcite{KB}. On se donne pour chaque paire une grande famille de chemins (des lignes brisées parallèles aux axes, avec un grand nombre de virages), on tire au hasard dans chaque famille, et quand des intersections se produisent, on les corrige localement.

Lorsque $m<n$ et $m$ impair, L. Guth combine un ingrédient géométrique, l'inégalité isopérimétrique, et un ingrédient topologique, le théorème de J.-P.~\textcite{Serre} sur la finitude des groupes d'homotopie des sphères : lorsque $m$ est impair, seul $\pi_m(S^m)$ est infini. Partant d'une homotopie, i.e. d'une application continue $h$ de la boule $B^{n+1}$ vers $S^m$, on choisit une triangulation de maille $1/L$, et on améliore $h$ squelette par squelette, de sorte que la restriction de $h$ à tout simplexe appartienne à une liste finie. Les obstructions rencontrées s'interprètent comme des cocycles, qu'on peut approcher quantitativement par des cobords bornés, en vertu de l'inégalité isopérimétrique, s'ils sont exacts. Quand ils ne le sont pas, c'est qu'on rencontre un élément non trivial d'un groupe d'homotopie de $S^m$, il y a à nouveau un nombre fini de modèles, d'après J.-P. Serre. En fin de récurrence, $h$ a une liste finie de modèle locaux, sa constante de Lipschitz est de l'ordre du diamètre de la boule triangulée, soit $L$.

Le théorème \ref{th:CDMW} s'applique aux H-espaces, et notamment aux classifiants. Il permet d'estimer la complexité des cobordismes. En effet, un théorème de R.~\textcite{Thom} ramène le calcul du groupe de cobordisme des variétés aux groupes d'homotopie de l'espace de Thom d'une grassmannienne. La complexité d'une variété (resp. d'une sous-variété) est mesurée ici par le volume $V$ d'une métrique riemannienne à géométrie locale bornée, dans le sens où le rayon d'injectivité est $\ge 1$, la courbure sectionnelle bornée par~$1$ (resp. en sus, rayon d'injectivité normal $\ge 1$, seconde forme fondamentale bornée par~$1$). Pour une variété à bord, on demande que les géométries locales de la variété et de son bord comme sous-variété soient bornées.

\begin{theo}[Appendice dans \textcite{CDMW}]\label{th:CDMWcob}
Il existe une fonction $\phi$ linéaire à un terme logarithmique près telle que, pour toute variété de complexité $V$ qui est le bord d'au moins une variété, borde aussi une variété de complexité majorée par $\phi(V)$.
\end{theo}

Cette borne presque linéaire est spectaculaire, elle donne envie d'utiliser le résultat pour des applications pratiques. 

La première étape de la démarche de R. Thom consiste à plonger la variété donnée, de dimension $n$, dans une sphère de dimension $N\ge 2n+1$. La version quantitative garantit la réalisation d'une variété de complexité $V$ comme sous-variété à géométrie locale bornée dans une boule euclidienne de rayon
$$
R\le C\,V^{1/(N-n)}(\log V)^{2n+2},
$$
soit très près de la borne inférieure donnée par des considérations de volume lorsque $N$~est grand, ou des considérations d'expansion pour la version simpliciale. Suivant \textcite{GG}, ce plongement ---disons, une version simpliciale--- est obtenu en deux temps. On tire d'abord au hasard les $V$ sommets sur la sphère de rayon $R_0=V^{1/(N-n)}$, parmi ceux dont les distances mutuelles sont $\ge \alpha R_0$ et tels que tous les simplexes de dimension $\le n$ construits sur l'échantillon sont $\alpha$-épais (pas de petits angles), pour un $\alpha>0$. Alors la probabilité qu'une boule unité intersecte un simplexe choisi à l'avance est $\le 1/V$. Comme il y a suffisamment d'indépendance partielle, on en déduit qu'avec forte probabilité, une boule unité ne rencontre qu'un nombre $O(\log V)$ de simplexes. Dans un second temps, on subdivise les simplexes de sorte que toutes les arêtes soient de longueur~$1$ et on améliore la géométrie locale, mais seulement à l'échelle $(\log V)^{-(2n+2)}$, faute de mieux. L'idée est que si l'on déplace au hasard d'une distance $\le 1$ les sommets de deux simplexes, alors leurs $\epsilon$-voisinages tubulaires deviennent disjoints avec probabilité au moins $1-C\epsilon$. L'exposant $2n+2$, c'est le nombre maximum de sommets cumulé de deux simplexes.

Noter que dès qu'il s'agit de plonger des objets $X$ de dimension $n$ dans $\R^{2n}$, les choses se gâtent. Pour certains complexes simpliciaux $X$, même si un plongement existe, sa complexité peut être exponentielle en la complexité de $X$, \parencite{FrKr}.

\subsection{Norme et familles de variétés}

Dès qu'on cherche à quantifier l'homotopie sur des familles d'espaces, la situation change, et des questions de complexité algorithmique ressurgissent, même en basse dimension. Soit $n\ge 3$ fixé. Soit $Y=S^n$ standard. On considère la famille~$\mathcal{X}$ des espaces~$X$ obtenus en triangulant~$S^n$. On mesure la taille d'un élément de~$\mathcal{X}$ par son nombre~$N$ de simplexes. On s'intéresse à la plus petite constante de Lipschitz d'une application de degré non nul $X\to Y$, notée $\nu{\not=0}(X)$.

\begin{theo}[\cite{BrGM}]\label{BrGM}
Calculer une approximation à $N^{\log\log(N)}$ près de la fonction $\nu{\not=0}$ sur la famille $\mathcal{X}$ est NP-difficile. 
\end{theo}

La démonstration consiste à construire une réduction polynômiale en $N$ du problème de la détermination des vecteurs les plus courts dans un réseau, dont la difficulté d'approximation a été démontrée par I.~\textcite{Dinur}. La réduction s'obtient en réalisant la sphère triangulée comme une sorte d'approximation métrique d'un polyèdre dont le $H^n$ est grand. 

En revanche, pour $n=2$, c'est-à-dire pour les $2$-sphères triangulées, un algorithme polynômial de calcul de la fonction $\nu{\not=0}$ est donné dans \textcite{Guth2005}.

\section{Conclusion}

Soulever des questions quantitatives en topologie, c'est ouvrir une boîte de Pandore : les questions fusent, mais elles sont difficiles. 

L'étude quantitative des homotopies a en particulier connu une très longue gestation. On a l'impression de n'avoir qu'effleuré la description du paysage de Morse de la fonction $\lip$ (on n'a pas dépassé les questions de connexité) et déjà les techniques remontant aux années 1970 ont dû être revisitées et adaptées. Du point de vue de la géométrie riemannienne, le rôle joué par les formes différentielles, objets pourtant linéaires, est frappant. Jusqu'où les formes différentielles peuvent-elles nous mener ? 

Ces évolutions plus récentes de la topologie que sont les géométries symplectique et de contact ont leur lot de questions quantitatives, très peu abordées jusqu'à présent.

Même si elles n'ont été abordées dans ce texte que de façon superficielle, les questions liées au calcul sont importantes. On a vu que la frontière entre décidabilité et indécidabilité se faufilait même en basses dimensions. On a rencontré presque tous les degrés de complexité algorithmique imaginables. De façon plus terre à terre, des outils de calcul de groupes d'homotopie, par exemple, sont disponibles mais encore à développer pour devenir de véritables assistants de preuve et entrer dans la boîte à outil quotidienne de la topologie. Les succès de l'analyse topologique des données reposent en partie sur sa capacité à contourner par des progrès conceptuels l'obstacle du coût du simple calcul de l'homologie à coefficients dans $\Z/2\Z$ de gros complexes simpliciaux, voir \textcite{Ginot}. 

On a du pain sur la planche.



\printbibliography

\end{document}


\textcite, \parencite

Courbe[(1 + (1/2)*sin(5*t))*cos(3*sin(t/2)-3*sin(t)), (1 + (1/2)*sin(5*t))*sin(3*sin(t/2)-3*sin(t)), t, 0, 2*Pi]

Courbe[(1 + 1 / 10 sin(4t)-(1/20)sin(t)) cos(2.5sin(t / 2) - 5sin(t)), (1 + 1 / 10 sin(4t)-(1/20)sin(t)) sin(2.5sin(t / 2) - 5sin(t)), t, 0, 6.28319]
t11=(9.5*6.28319/10)
A_{11}=((1 + 1 / 10 sin(4t11)-(1/20)sin(t11)) cos(2.5sin(t11 / 2) - 5sin(t11)), (1 + 1 / 10 sin(4t11)-(1/20)sin(t11)) sin(2.5sin(t11 / 2) - 5sin(t11)))

https://www.youtube.com/watch?v=Vp_NrF9zfEw
7:50 à 8:26
